\documentclass{article}%
\usepackage{amsmath}
\usepackage{amsfonts}
\usepackage{amssymb}
\usepackage{graphicx}%
\providecommand{\U}[1]{\protect \rule{.1in}{.1in}}
\newtheorem{theorem}{Theorem}

\newtheorem{lemma}[theorem]{Lemma}

\newtheorem{remark}[theorem]{Remark}

\begin{document}


\title{Explicit solutions of $G$-heat equation with a class of initial conditions by $G$-Brownian motion}
\author{Mingshang Hu\footnote{Email: humingshang2004\symbol{64}hotmail.com}\\
{\small School of Mathematics, Shandong University, Jinan,
\small{250100}}}
\date{\today}
\maketitle

\noindent\textbf{Abstract. }
We obtain the viscosity solution of $G$-heat equation with the
initial condition $\phi(x)=x^{n}$ for each integer $n\geq1$ using
the method of $G$-Brownian motion.

\  \  \

\noindent\textbf{Keywords: }$G$-heat equation, sublinear
expectation, $G$-normal distribution, $G$-expectation, $G$-Brownian
motion.

\  \  \

\section{Introduction}

The notions of $G$-normal distribution, $G$-expectation and
$G$-Brownian motion were firstly introduced by Peng (see
\cite{Peng-3} and \cite{Peng-4}) via the following $G$-heat
equation:
\begin{equation}
\label{e1} \frac{\partial u}{\partial
t}-G(\frac{\partial^{2}u}{\partial x^{2}})=0, \quad u|_{t=0}=\phi,
\quad (t,x)\in [0, \infty)\times\mathbb{R},
\end{equation}
where $G(a)=\frac{1}{2}(a^{+}-\sigma^{2}a^{-})$. Here $\sigma\in
[0,1]$ is a fixed constant and $a^{+}=\max\{0, a \}$,
$a^{-}=(-a)^{+}$. This equation is also the Barenblatt equation
except the case $\sigma=0$ (see \cite{Ba}, \cite{BS} and
\cite{KPV}).

Under the sublinear framework, Peng gave many important notions
corresponding to linear case and obtained many important properties
of $G$-Brownian motion (see \cite{Peng-6}). Recently, Peng developed
the law of large numbers and central limit theorem under sublinear
expectations, which indicate that $G$-normal distribution plays the
same important role in the theory of sublinear expectations as
normal distribution in the linear expectations (see \cite{Peng-5}
and \cite{Peng-7}).

Since the importance of $G$-normal distribution and $G$-Brownian
motion, Peng proposed the problem of how to calculate the
$G$-expectation of $\phi(B_t)$? For convex or concave $\phi$, Peng
gave the formula in \cite{Peng-3} and \cite{Peng-4}. But for neither
convex nor concave $\phi$, how to calculate? In particular, the
calculation of $B_t^{2n+1}$, for each integer $n\geq 1$.

In this paper, we give the relation between the solution of $G$-heat
equation (\ref{e1}) with the initial condition $\phi(x)=x^{n}$ for
each integer $n\geq1$ and the solution of ordinary differential
equation (see Section 3), then we can get the solution of $G$-heat
equation (\ref{e1}) by solving ordinary differential equation. In
particular, we get the $G$-expectation of $B_t^{2n+1}$ for each
integer $n\geq 1$. In fact, we also get the solution of the
Barenblatt equation with the same initial condition. Finally, we
point out the application of our result in mathematical finance.

This paper is organized as follows: in Section 2, we recall briefly
the notions of $G$-normal distribution, $G$-expectation and
$G$-Brownian motion. Also, we prove that the $G$-expectation of
$B_t^{2n+1}$ is strictly bigger than 0 which is different from
linear case. The main theorem in which we get the solution of
$G$-heat equation (\ref{e1}) with the initial condition
$\phi(x)=x^{n}$ for each integer $n\geq1$ is stated and proved in
Section 3. In Section 4, we consider the application of our result
in mathematical finance.

\section{$G$-Brownian Motions under $G$-expectations}

Let $\Omega=C_0(\mathbb{R}^{+})$ be the space of all
$\mathbb{R}$-valued continuous paths
$(\omega_t)_{t\in\mathbb{R}^+}$, with $\omega_0=0$, equipped with
the distance
\[
\rho(\omega^1,\omega^2):=\sum_{i=1}^{\infty}2^{-i}[(\max_{t\in
[0,i]}|\omega_t^1-\omega_t^2|)\wedge1].
\]
The corresponding canonical process $B_t(\omega)=\omega_t$, $t\in
[0,\infty)$, for $\omega\in \Omega$. We denote
$C_{l,Lip}(\mathbb{R}^n)$ the linear space of (local Lipschitz)
functions $\phi$ satisfying
\begin{align*}
|\phi(x)-\phi(y)| &\leq C(1+|x|^m+|y|^m)|x-y|,\quad \forall x, y\in
\mathbb{R}^n,\\
& \mbox{for some}\  C >0, m\in \mathbb{N} \ \mbox{depending on} \
\phi.
\end{align*}
For each fixed $T\in [0,\infty)$, we set
\[
L_{ip}^0(\mathcal{F}_T):=\{\phi(B_{t_1},\cdots,B_{t_n}):\forall\,n\in\mathbb{N},\,t_1,\cdots,t_n\in
[0,T],\,\forall\,\phi\in C_{l,Lip}(\mathbb{R}^n) \}.
\]
It is clear that $L_{ip}^0(\mathcal{F}_t)\subseteq
L_{ip}^0(\mathcal{F}_T)$, for $t\leq T$. We also set
\[
L_{ip}^0(\mathcal{F}):=\bigcup_{n=1}^{\infty}L_{ip}^0(\mathcal{F}_n).
\]
Following \cite{Peng-3} and \cite{Peng-4}, we can construct a
consistent sublinear expectation called $G$-expectation
$\hat{\mathbb{E}}[\cdot]: L_{ip}^0(\mathcal{F})\mapsto \mathbb{R}$
satisfying, for each $X,Y \in L_{ip}^0(\mathcal{F})$,
\begin{description}
  \item[(a)]\textbf{Monotonicity}: if $X\geq Y$, then
  $\hat{\mathbb{E}}[X]\geq\hat{\mathbb{E}}[Y]$.
  \item[(b)]\textbf{Constant preserving}: $\hat{\mathbb{E}}[c]=c$, $\forall\, c\in\mathbb{R}$.
  \item[(c)]\textbf{Sub-additivity}:
  $\hat{\mathbb{E}}[X+Y]\leq\hat{\mathbb{E}}[X]+\hat{\mathbb{E}}[Y]$.
  \item[(d)]\textbf{Positive homogeneity}: $\hat{\mathbb{E}}[\lambda
  X]=\lambda\hat{\mathbb{E}}[X]$, $\forall\,\lambda\geq0$.
\end{description}
Under $G$-expectation $\hat{\mathbb{E}}[\cdot]$, the canonical
process $\{ B_t:t\geq 0 \}$ is called $G$-Brownian motion, and the
distribution of $B_1$ is said to be the $G$-normal distribution,
i.e., for each $\phi\in C_{l,Lip}(\mathbb{R})$, the function
\[
u(t,x):=\hat{\mathbb{E}}[\phi(x+\sqrt{t}B_1)],\quad (t,x)\in
[0,\infty)\times\mathbb{R}
\]
is the viscosity solution (see \cite{CIL}) of $G$-heat equation
(\ref{e1}).
Moreover, $G$-Brownian motion has independent and stationary
increments, i.e., for each $n\in \mathbb{N}$, $\phi\in
C_{l,Lip}(\mathbb{R}^n)$, $s,t \geq0$ and $t_1,\ldots,t_{n-1}\in
[0,t]$, we have
\begin{equation} \label{e2}
\hat{\mathbb{E}}[\phi(B_{t_1},\ldots,B_{t_{n-1}},B_{t+s}-B_{t})]=\hat{\mathbb{E}}[\varphi(B_{t_1},\ldots,B_{t_{n-1}})],
\end{equation}
where
$\varphi(x_1,\ldots,x_{n-1})=\hat{\mathbb{E}}[\phi(x_1,\ldots,x_{n-1},\sqrt{s}B_1)]$.
Specially, for each $\phi\in C_{l,Lip}(\mathbb{R})$,
$\hat{\mathbb{E}}[\phi(B_t)]=\hat{\mathbb{E}}[\phi(\sqrt{t}B_1)]$,
$\forall\,t>0$. In particular, for each integer $n\geq1$,
$\hat{\mathbb{E}}[B_t^{2n+1}]=t^{n+\frac{1}{2}}\hat{\mathbb{E}}[B_1^{2n+1}]$.

\begin{remark}
For $\sigma\in (0,1]$, $G$-heat equation (\ref{e1}) is a uniformly
parabolic PDE and $G$ is a convex function, then it has a unique
$C^{1,2}$ solution (see \cite {Kr} and \cite {Wang}).
\end{remark}

In the following, $P^W$ denotes Wiener probability measure on
$\Omega$, $E^W$ always denotes the linear expectation with respect
to $P^W$. Under $P^W$, the canonical process $\{ B_t:t\geq 0 \}$ is
the classical standard Brownian motion.

\begin{lemma}
For each fixed $\sigma\in [0,1)$, we have
\begin{description}
  \item[(i)] For each $\phi\in C_{l,Lip}(\mathbb{R})$, $\hat{\mathbb{E}}[\phi(B_t)]\geq \sup_{\sigma\leq \nu
\leq1}E^W[\phi(\nu B_t)]$.
  \item[(ii)] For each integer $n\geq1$, $\hat{\mathbb{E}}[B_1^{2n+1}]>0$.
\end{description}
\end{lemma}

\noindent\textbf{Proof. } Noting that
$u(t,x):=\hat{\mathbb{E}}[\phi(x+B_t)]$ and
$u^\nu(t,x):=E^W[\phi(x+\nu B_t)]$ are respective the viscosity
solution of $G$-heat equation (\ref{e1}) and
$u^\nu_t-\frac{1}{2}\nu^2u^\nu_{xx}=0$ with the same initial
condition, by comparison theorem for parabolic partial differential
equations (see \cite{CIL}), we get (i). It follows from (\ref{e2})
that
\[
\hat{\mathbb{E}}[B_1^{2n+1}]=\hat{\mathbb{E}}[(B_{\frac{1}{2}}+B_1-B_{\frac{1}{2}})^{2n+1}]
=\hat{\mathbb{E}}[\phi(B_{\frac{1}{2}})],
\]
where $\phi(x)=\hat{\mathbb{E}}[(x+B_{\frac{1}{2}})^{2n+1}]$. As a
consequence of (i),
\begin{align*}
\phi(x) & \geq \sup_{\sigma\leq \nu\leq 1}E^W[(x+\nu
B_{\frac{1}{2}})^{2n+1}]\\
 & \geq
\sum_{i=0}^n\binom{2n+1}{2i}x^{2(n-i)+1}E^W[B_{\frac{1}{2}}^{2i}]
+\frac{n(2n+1)}{2}(1-\sigma^2)(x^-)^{2n-1}.
\end{align*}
Therefore,
\[
\hat{\mathbb{E}}[B_1^{2n+1}]=\hat{\mathbb{E}}[\phi(B_{\frac{1}{2}})]\geq
E^W[\phi(B_{\frac{1}{2}})]\geq
\frac{n(2n+1)}{2}(1-\sigma^2)E^W[(B_{\frac{1}{2}}^-)^{2n-1}].
\]
As $\sigma<1$, we get $\hat{\mathbb{E}}[B_1^{2n+1}]>0$ for each
integer $n\geq1$. $\Box$

\begin{remark}
For convex (resp. concave) $\phi\in C_{l,Lip}(\mathbb{R})$, we have
$\hat{\mathbb{E}}[\phi(B_t)]=E^W[\phi(B_t)]$ (resp.
$\hat{\mathbb{E}}[\phi(B_t)]=E^W[\phi(\sigma B_t)]$); for odd
function $\phi\in C_{l,Lip}(\mathbb{R})$, we have
$\hat{\mathbb{E}}[\phi(B_t)]=\hat{\mathbb{E}}[-\phi(B_t)]$, all see
Proposition 7 in \cite{Peng-3}. For $\sigma<1$, $t>0$,
$\hat{\mathbb{E}}[B_t^{2n+1}]=\hat{\mathbb{E}}[-B_t^{2n+1}]>0$, this
is different from linear case, and at the same time show that the
inequality in (i) can't be changed into equality.
\end{remark}

\begin{remark}
We define $\mathcal{F}_{t}:=\sigma \{B_{u},0\leq u\leq t\} \vee
\mathcal{N}^W$ for each $t\geq 0$, where $\mathcal{N}^W$ is all
$P^W$-null sets, $\mathcal{V}$ is the set of all
$\mathcal{F}_{t}$-adapted controls with values in $[\sigma,1]$.
Using a stochastic control method (see \cite{DHP}), we have for each
$\phi\in C_{l,Lip}(\mathbb{R})$,
\[
\hat{\mathbb{E}}[\phi(B_t)]=\sup_{\nu_. \in
\mathcal{V}}E^W[\phi(\int_0^t\nu_s\,dB_s)].
\]
\end{remark}

\section{Main Result}

In this Section, we discuss the solution of $G$-heat equation
(\ref{e1}) with the initial condition $\phi(x)=x^{n}$ for each
integer $n\geq1$. Noting that for convex (resp. concave) $\phi$,
$u(t,x):=E^W[\phi(x+B_t)]$ (resp. $u(t,x):=E^W[\phi(x+\sigma B_t)]$)
is the solution of $G$-heat equation (\ref{e1}). Consequently, we
only consider the following $G$-heat equation:
\begin{equation} \label{e3}
\frac{\partial u}{\partial t}-\frac{1}{2}\big(
(\frac{\partial^{2}u}{\partial
x^{2}})^{+}-\sigma^2(\frac{\partial^{2}u}{\partial x^{2}})^{-} \big
)=0, \quad u(0,x)=x^{2n+1},
\end{equation}
where $\sigma\in [0,1]$, $(t,x)\in [0, \infty)\times\mathbb{R}$. For
$\sigma=1$, $G$-heat equation (\ref{e3}) is classical heat equation,
then $u(t,x):=E^W[(x+B_t)^{2n+1}]$ is the unique solution. In order
to give the solution of $G$-heat equation (\ref{e3}) with $\sigma\in
[0,1)$, we define for each integer $n\geq0$,
\begin{equation} \label{e4}
\begin{array}{r@{\;=\;}l}
\displaystyle g_n(x)&
\displaystyle\sum_{i=0}^{n}\frac{(2n+1)!}{(2(n-i))!!(2i+1)!}x^{2i+1},
\\[0.5cm]
\displaystyle h_n(x)&\displaystyle
\sum_{i=0}^{n}\frac{(n+i)!(n-i)!}{(2(n-i))!!(2i)!}\left(
\binom{2n+1}{0}+\cdots+\binom{2n+1}{n-i} \right)x^{2i} ,
\end{array}
\end{equation}
with the convention that $0!=1$ and $0!!=1$. Now, we give our main
Theorem.

\begin{theorem}
For each integer $n\geq1$, $g_n(\cdot)$ and $h_n(\cdot)$ defined in
(\ref{e4}). The following statements hold.\\
{\bf (i)} For each fixed $\sigma\in (0,1)$, we define
\begin{equation}\label{e5}
P_n^{\sigma}(x)\!=\!\!\left \{\!\!\!
\begin{array}{l}
g_n(x)+\frac{k_n}{(2n)!!}[h_n(x)\exp(-\frac{x^2}{2})
-g_n(x)\int_x^{\infty}
\exp(-\frac{t^2}{2})\,dt \,],\qquad\qquad x\geq c_n\\[2.3mm]

\sigma^{2n+1}g_n(\frac{x}{\sigma})+\frac{d_n}{(2n)!!}[h_{n}(\frac{x}{\sigma})
\exp(-\frac{x^2}{2\sigma^2})+g_{n}(\frac{x}{\sigma})
\int_{-\infty}^{\frac{x}{\sigma}}\exp(-\frac{t^2}{2})\,dt \,],x<c_n
\end{array}
 \right.
\end{equation}
where $c_n$, $k_n$ and $d_n$ are constants such that
\begin{equation}\label{e6}
\left \{
\begin{array}{l}
\displaystyle
h_{n-1}(\frac{c_n}{\sigma})+g_{n-1}(\frac{c_n}{\sigma})\exp(\frac{c_n^2}{2\sigma^2})
\int_{-\infty}^{\textstyle\frac{c_n}{\sigma}}\exp(-\frac{t^2}{2})\,dt
\\[4mm]
\displaystyle \qquad \qquad \qquad =\sigma^{2n}[\,h_{n-1}(c_n)
-g_{n-1}(c_n)\exp(\frac{c_n^2}{2})\int_{c_n}^{\infty}\exp(-\frac{t^2}{2})\,dt\,] \\[4mm]
\displaystyle
k_n=-\frac{(2n)!!g_{n-1}(c_n)}{h_{n-1}(c_n)\exp(-\frac{c_n^2}{2})
-g_{n-1}(c_n)\int_{c_n}^{\infty}\exp(-\frac{t^2}{2})\,dt}\\[4mm]
\displaystyle d_n=\frac{\sigma
g_{n-1}(\frac{c_n}{\sigma})}{g_{n-1}(c_n)}
k_n\exp(\frac{1-\sigma^2}{2\sigma^2}c_n^2)\\[4mm]
\displaystyle c_n<0
\end{array}
\right.
\end{equation}
Then $\displaystyle
u_n^{\sigma}(t,x):=t^{n+\frac{1}{2}}P_n^{\sigma}(\frac{x}{\sqrt{t}})$
is the unique solution of $G$-heat equation (\ref{e3}). \\
{\bf (ii)} For $\sigma=0$, we define
\begin{equation}\label{e7}
P_n(x)=\left \{
\begin{array}{l}
\displaystyle
g_n(x)+\frac{k_n}{(2n)!!}[\,h_n(x)\exp(-\frac{x^2}{2})-g_n(x)\int_x^{\infty}
\exp(-\frac{t^2}{2})\,dt \,],x\geq c_n\\[2mm]
\displaystyle x^{2n+1}, \qquad \qquad \qquad \qquad \qquad \qquad
\qquad \qquad \qquad \qquad \qquad  x<c_n
\end{array}
 \right.
\end{equation}
where $c_n$, $k_n$ are constants such that
\begin{equation}\label{e8}
\left \{
\begin{array}{l}
\displaystyle (2n-1)!=c_n^{2n}[\,h_{n-1}(c_n)
-g_{n-1}(c_n)\exp(\frac{c_n^2}{2})\int_{c_n}^{\infty}\exp(-\frac{t^2}{2})\,dt\,]
\\[2mm]
\displaystyle
k_n=-\frac{2n}{(2n-1)!!}g_{n-1}(c_n)c_n^{2n}\exp(\frac{c_n^2}{2})\\[2mm]
\displaystyle c_n<0
\end{array}
\right.
\end{equation}
Then $\displaystyle
u_n(t,x):=t^{n+\frac{1}{2}}P_n(\frac{x}{\sqrt{t}})$ is the unique
viscosity solution of $G$-heat equation (\ref{e3}).
\end{theorem}

In order to prove the main Theorem, we need the following Lemmas.
The first Lemma gives the relation between the solution of $G$-heat
equation (\ref{e3}) and the solution of the following ordinary
differential equation (ODE for short):
\begin{equation}\label{e9}
(y'')^+-\sigma^2(y'')^-+xy'-(2n+1)y=0,\quad
\lim_{t\rightarrow0^+}t^{n+\frac{1}{2}}y(\frac{x}{\sqrt{t}})=x^{2n+1}.
\end{equation}

\begin{lemma}
Let $\sigma\in (0,1]$ be a fixed constant. The following statements
hold.
\begin{description}
  \item[(i)] If $u$ is the solution of $G$-heat equation (\ref{e3}), then
$P(x):=u(1,x)$ is the $C^{\:2}$ solution of ODE (\ref{e9}).
  \item[(ii)] If $P$ is the $C^{\:2}$ solution of ODE (\ref{e9}), then
$\displaystyle u(t,x):=t^{n+\frac{1}{2}}P(\frac{x}{\sqrt{t}})$ is
the solution of $G$-heat equation (\ref{e3}).
\end{description}
\end{lemma}

\noindent\textbf{Proof. }For (i), noting that $B_1$ is $G$-normal
distributed, then for $t>0$,
\[
u(t,x)=\hat{\mathbb{E}}[(x+\sqrt{t}B_1)^{2n+1}]
=t^{n+\frac{1}{2}}\hat{\mathbb{E}}[(\frac{x}{\sqrt{t}}+B_1)^{2n+1}]
=t^{n+\frac{1}{2}}u(1,\frac{x}{\sqrt{t}}).
\]
Consequently, $\displaystyle
u(t,x)=t^{n+\frac{1}{2}}P(\frac{x}{\sqrt{t}})$ for $t>0$. Also
$u(0,x)=x^{2n+1}$, we get
$\lim_{t\rightarrow0^+}t^{n+\frac{1}{2}}P(\frac{x}{\sqrt{t}})=x^{2n+1}$.
According to Remark $1$, we have $P\in C^{\:2}$. It is easy to
verify that, for $t>0$,
\begin{equation}\label{e10}
\frac{\partial u}{\partial t}=t^{n-\frac{1}{2}} [\
(n+\frac{1}{2})P(\frac{x}{\sqrt{t}})-\frac{x}{2\sqrt{t}}P'(\frac{x}{\sqrt{t}})\
], \quad \frac{\partial^{2}u}{\partial
x^{2}}=t^{n-\frac{1}{2}}P''(\frac{x}{\sqrt{t}}).
\end{equation}
Substituting (\ref{e10}) into $G$-heat equation (\ref{e3}), we
conclude the result. The proof for (ii) is similar. $\Box$

\begin{remark}
Let $u^{\sigma}$ denote the viscosity solution of $G$-heat equation
(\ref{e3}) with $\sigma\in [0,1]$, then for each fixed $(t,x)\in (0,
\infty)\times\mathbb{R} $, $\phi(\sigma):=u^{\sigma}(t,x), \sigma\in
[0,1]$, is a continuous and decreasing function (see \cite{CIL} and
\cite{Peng-3}).
\end{remark}

In the following, we first solve ODE (\ref{e9}) for $\sigma\in
(0,1]$, then by Lemma $(6)$ we get the corresponding solution of
$G$-heat equation (\ref{e3}). By the above Remark, taking
$\sigma\downarrow 0$, we get the solution of $G$-heat equation
(\ref{e3}) with $\sigma=0$. For this, we give the following Lemmas.
%

\begin{lemma}
Let $\sigma\not=0$, $\alpha>0$, if $\varphi(\cdot)$ is a solution of
ODE $y''+xy'-\alpha y=0$, then $\psi(x):=\varphi(\frac{x}{\sigma})$
is a solution of ODE $\sigma^2y''+xy'-\alpha y=0$.
\end{lemma}

\textbf{Proof. }It is easy to verify the result. $\Box$

\  \  \

In the following, we use $g_n(\cdot)$ and $h_n(\cdot)$
defined in (\ref{e4}).

\begin{lemma}
The general solution of ODE $y''+xy'-(2n+1)y=0$ is
\[
y(x)=\lambda_1g_n(x)+\lambda_2[\,h_n(x)\exp(-\frac{x^2}{2})-g_n(x)\int_x^{\infty}
\exp(-\frac{t^2}{2})\,dt \,],
\]
where $\lambda_1$ and $\lambda_2$ are arbitrary constants.
\end{lemma}

\noindent\textbf{Proof. }Applying Lemma $6$ with $\sigma=1$, we get
$\bar{y}_n(x)=E^W[(x+B_1)^{2n+1}]$ and
$\hat{y}_n(x)=\sqrt{2\pi}E^W[((x+B_1)^-)^{2n+1}]$ are two linear
independent solutions of ODE $y''+xy'-(2n+1)y=0$. It is easy to
check that $\bar{y}_n(x)=g_n(x)$ and
\[
\hat{y}_n(x)=h_n(x)\exp(-\frac{x^2}{2})-g_n(x)\int_x^{\infty}
\exp(-\frac{t^2}{2})\,dt,
\]
which completes the proof. $\Box$

\begin{lemma}
Let $x>0$, then for each integer $n\geq1$, we have
\begin{equation}\label{e11}
\frac{h'_{n-1}(x)+g_{n-1}(x)}{g'_{n-1}(x)+xg_{n-1}(x)}\exp(-\frac{x^2}{2})\leq
\int_x^{\infty}\exp(-\frac{t^2}{2})\,dt\leq
\frac{h_n(x)}{g_n(x)}\exp(-\frac{x^2}{2}).
\end{equation}
\end{lemma}

\noindent\textbf{Proof. }We define
$m_n(x)=\sqrt{2\pi}E^W[((x+B_1)^-)^{2n+1}]$, then by the proof of
the above Lemma, we know
\[
m_n(x)=h_n(x)\exp(-\frac{x^2}{2})-g_n(x)\int_x^{\infty}
\exp(-\frac{t^2}{2})\,dt.
\]
Noting $m_n(x)\geq 0$, the right-hand side inequality of (\ref{e11})
holds. We also observe that $m_n(\cdot)$ is a decreasing function,
then $m_n'(x)\leq 0$, which yields the left-hand side inequality of
(\ref{e11}). $\Box$

\begin{remark}
It is easy to verify that, for each fixed $x>0$, as
$n\rightarrow\infty$,
\[
\frac{h'_{n-1}(x)+g_{n-1}(x)}{g'_{n-1}(x)+xg_{n-1}(x)}\exp(-\frac{x^2}{2})\Big\uparrow
\int_x^{\infty}\exp(-\frac{t^2}{2})\,dt,
\]
\[
\frac{h_n(x)}{g_n(x)}\exp(-\frac{x^2}{2})\Big\downarrow
\int_x^{\infty}\exp(-\frac{t^2}{2})\,dt.
\]
Furthermore, we can also prove that
\[
\int_x^{\infty}\exp(-\frac{t^2}{2})\,dt-\frac{h'_{n-1}(x)+g_{n-1}(x)}{g'_{n-1}(x)+xg_{n-1}(x)}\exp(-\frac{x^2}{2})<\min_{1\leq
k\leq n}\{\frac{\sqrt{2\pi}(2k)!(n-k)!}{x^{2k}2^{3k}n!}  \},
\]
\[
\frac{h_n(x)}{g_n(x)}\exp(-\frac{x^2}{2})-\int_x^{\infty}\exp(-\frac{t^2}{2})\,dt<\min_{1\leq
k\leq
n-1}\{\frac{\sqrt{2\pi}(2k+1)!(2(n-k)-1)!n!}{2^kx^{2(k+1)}(2n)!(n-k-1)!}
\}.
\]
\end{remark}

\begin{lemma}
Let $\sigma\in (0,1)$ be a fixed constant, and we define for
$n\geq1$,
\begin{align*}
f_n(x):=&h_{n-1}(\frac{x}{\sigma})+g_{n-1}(\frac{x}{\sigma})\exp(\frac{x^2}{2\sigma^2})
\int_{-\infty}^{\textstyle\frac{x}{\sigma}}\exp(-\frac{t^2}{2})\,dt-\sigma^{2n}h_{n-1}(x)\\
&+\sigma^{2n}g_{n-1}(x)\exp(\frac{x^2}{2})\int_x^{\infty}\exp(-\frac{t^2}{2})\,dt,
\end{align*}
then there exists a unique $x_0<0$ such that $f_n(x_0)=0$.
\end{lemma}

\noindent\textbf{Proof. }It is easy to verify that
\[
f_n(0)=(1-\sigma^{2n})(2(n-1))!!>0,\quad f_n(-\infty)=-\infty,
\]
\begin{align*}
f_n'(x)=&\frac{1}{\sigma}(g_{n-1}'(\frac{x}{\sigma})+\frac{x}{\sigma}g_{n-1}(\frac{x}{\sigma}))\exp(\frac{x^2}{2\sigma^2})
\int_{-\infty}^{\textstyle\frac{x}{\sigma}}\exp(-\frac{t^2}{2})\,dt\\
&+\frac{1}{\sigma}(h_{n-1}'(\frac{x}{\sigma})+g_{n-1}(\frac{x}{\sigma}))-\sigma^{2n}(h_{n-1}'(x)+g_{n-1}(x)) \\
&+\sigma^{2n}(g_{n-1}'(x)+xg_{n-1}(x))\exp(\frac{x^2}{2})\int_x^{\infty}\exp(-\frac{t^2}{2})\,dt.
\end{align*}
As $x<0$,
\begin{align*}
f_n'(x)>&\frac{1}{\sigma}(g_{n-1}'(\frac{x}{\sigma})+\frac{x}{\sigma}g_{n-1}(\frac{x}{\sigma}))\exp(\frac{x^2}{2\sigma^2})
\int_{-\infty}^{\textstyle\frac{x}{\sigma}}\exp(-\frac{t^2}{2})\,dt\\
&+\frac{1}{\sigma}(h_{n-1}'(\frac{x}{\sigma})+g_{n-1}(\frac{x}{\sigma})).
\end{align*}
Applying Lemma 10, we get $f_n'(x)>0$ for $x<0$, which completes the
proof. $\Box$

\begin{lemma}
For each integer $n\geq1$, we have
\[
\lim_{x\rightarrow -\infty}x^{2n}[\, h_{n-1}(x)
+g_{n-1}(x)\exp(\frac{x^2}{2})\int_{-\infty}^{x}\exp(-\frac{t^2}{2})\,dt\,]=(2n-1)!.
\]
\end{lemma}

\noindent\textbf{Proof. }It is easy to prove that
\[
h_{n-1}(x)g_n(x)-g_{n-1}(x)h_n(x)=(2n-1)!x,
\]
\[
h_{n-1}(x)g_{n-1}'(x)+xh_{n-1}(x)g_{n-1}(x)-h_{n-1}'(x)g_{n-1}(x)-g_{n-1}^2(x)=(2n-1)!,
\]
then for $x<0$, by Lemma 10, we obtain
\begin{align*}
x^{2n}[\, h_{n-1}(x)
+g_{n-1}(x)\exp(\frac{x^2}{2})\int_{-\infty}^{x}\exp(-\frac{t^2}{2})\,dt\,]&\geq
\frac{(2n-1)!x^{2n+1}}{g_n(x)},\\
x^{2n}[\, h_{n-1}(x)
+g_{n-1}(x)\exp(\frac{x^2}{2})\int_{-\infty}^{x}\exp(-\frac{t^2}{2})\,dt\,]
& \leq\frac{(2n-1)!x^{2n}}{g_{n-1}'(x)+xg_{n-1}(x)}.
\end{align*}
Taking $x\rightarrow -\infty$ yields the result. $\Box$

\  \   \

\noindent\textbf{Proof of Theorem 5. }For each fixed $\sigma\in
(0,1)$, we denote by $P_n^{\sigma}(\cdot)$ the $C^{\:2}$ solution of
ODE (\ref{e9}). We assume that there exists a constant $c_n<0$ such
that
\[
(\frac{\mathrm{d}^2}{\mathrm{d}x^2}P_n^{\sigma})(x)\geq 0\
\mbox{for}\ x\geq c_n;\quad
(\frac{\mathrm{d}^2}{\mathrm{d}x^2}P_n^{\sigma})(x)\leq 0\
\mbox{for}\ x\leq c_n.
\]
Under this assumption, by Lemma 8, Lemma 9 and
$\lim_{t\rightarrow0^+}t^{n+\frac{1}{2}}P_n^{\sigma}(\frac{x}{\sqrt{t}})=x^{2n+1}$,
it follows that $P_n^{\sigma}(\cdot)$ has an expression as in
(\ref{e5}), where $c_n$, $k_n$ and $d_n$ are undetermined constants.
Note that $P_n^{\sigma}\in C^{\:2}$, we have
\[
\lim_{x\rightarrow c_n^{+}}P_n^{\sigma}(x)=\lim_{x\rightarrow
c_n^{-}}P_n^{\sigma}(x),\quad \lim_{x\rightarrow
c_n^{+}}(\frac{\mathrm{d}^2}{\mathrm{d}x^2}P_n^{\sigma})(x)=0,\quad
\lim_{x\rightarrow
c_n^{-}}(\frac{\mathrm{d}^2}{\mathrm{d}x^2}P_n^{\sigma})(x)=0.
\]
This yields (\ref{e6}) by direct verification. Applying Lemma 12, we
know that there exists a unique $c_n<0$ satisfying the first
equation in (\ref{e6}), then $k_n$ and $d_n$ are also unique. It is
easy to verify that $P_n^{\sigma}$ determined by (\ref{e5}) and
(\ref{e6}) belongs to $C^{\:2}$. Now, we show that this
$P_n^{\sigma}$ satisfies the above assumption. For $x\geq c_n$,
\begin{align*}
&\frac{1}{(2n+1)(2n)}(\frac{\mathrm{d}^2}{\mathrm{d}x^2}P_n^{\sigma})(x)\\
&=g_{n-1}(x)+\frac{k_n}{(2n)!!}(h_{n-1}(x)\exp(-\frac{x^2}{2})-g_{n-1}(x)\int_x^{\infty}\exp(-\frac{t^2}{2})\,dt)\\
&=E^W[(x+B_1)^{2n-1}]-\frac{E^W[((x+B_1)^-)^{2n-1}]}{E^W[((c_n+B_1)^-)^{2n-1}]}
E^W[(c_n+B_1)^{2n-1}] \\
&=E^W[((x+B_1)^+)^{2n-1}]-\frac{E^W[((x+B_1)^-)^{2n-1}]}{E^W[((c_n+B_1)^-)^{2n-1}]}
E^W[((c_n+B_1)^+)^{2n-1}] .
\end{align*}
It is easy to see that, for $x\geq c_n$,
\[
E^W[((c_n+B_1)^-)^{2n-1}]\geq E^W[((x+B_1)^-)^{2n-1}],
\]
\[
E^W[((x+B_1)^+)^{2n-1}]\geq E^W[((c_n+B_1)^+)^{2n-1}],
\]
therefore, we get
$(\frac{\mathrm{d}^2}{\mathrm{d}x^2}P_n^{\sigma})(x)\geq 0$ for
$x\geq c_n$. Similarly, we can prove that
$(\frac{\mathrm{d}^2}{\mathrm{d}x^2}P_n^{\sigma})(x)\leq 0$ for
$x\leq c_n$. Thus, this $P_n^{\sigma}$ is indeed the $C^{\:2}$
solution of ODE (\ref{e9}). By Lemma 6,
$u_n^{\sigma}(t,x):=t^{n+\frac{1}{2}}P_n^{\sigma}(\frac{x}{\sqrt{t}})$
is the solution of $G$-heat equation (\ref{e3}), and then the proof
of part {\bf (i)} is complete. We now prove part {\bf (ii)}. We
define
\[
l_n(x)=\int_x^{\infty}\exp(-\frac{t^2}{2})\,dt-\frac{h_{n-1}(x)}{g_{n-1}(x)}\exp(-\frac{x^2}{2}),\quad
x<0.
\]
It is easy to prove that
\begin{equation}\label{e12}
\lim_{x\rightarrow0^-}l_n(x)=+\infty;\quad
l_n'(x)=\frac{(2n-1)!}{g_{n-1}^2(x)}\exp(-\frac{x^2}{2})>0,\ x<0.
\end{equation}
According to Remark 7,
$u^{\sigma}(1,0)=k_n(\sigma)=\frac{(2n)!!}{l_n(c_n(\sigma))}\downarrow0
$ as $\sigma\uparrow 1$, then by (\ref{e12}) we get
$c_n(\sigma)\uparrow 0$ as $\sigma\uparrow 1$. Hence,
$\frac{c_n(\sigma)}{\sigma}\downarrow -\infty$ as
$\sigma\downarrow0$. Taking $\sigma\downarrow 0$, by Lemma 13, we
obtain {\bf (ii)}. $\Box$

\begin{remark}
For each fixed $\sigma\in (0,1)$,
$\hat{\mathbb{E}}[B_t^{2n+1}]=u_n^{\sigma}(t,0)=k_n\,t^{n+\frac{1}{2}}$,
where $k_n$ is determined by (\ref{e6}). In order to get $k_n$, we
solve the first equation in (\ref{e6}) of $c_n$, then $k_n$ is
uniquely determined by $c_n$. Specially,
$\hat{\mathbb{E}}[B_t^{3}]=k\,t^{\frac{3}{2}}$, where $k$ is
determined by the following equations:
\begin{equation}
 \left \{
\begin{array}{l}
\displaystyle 1+\frac{c}{\sigma}\exp(\frac{c^2}{2\sigma^2})
\int_{-\infty}^{\textstyle\frac{c}{\sigma}}\exp(-\frac{t^2}{2})\,dt\\[4mm]
\displaystyle \qquad\qquad\qquad\qquad
=\sigma^2[\,1-c\exp(\frac{c^2}{2})\int_c^{\infty}\exp(-\frac{t^2}{2})\,dt\,]
\\[4mm]
\displaystyle
k=-\frac{2c}{\exp(-\frac{c^2}{2})-c\int_c^{\infty}\exp(-\frac{t^2}{2})\,dt}\\[4mm]
\displaystyle c<0
\end{array}
\right.
\end{equation}
For $\sigma=0$,
$\hat{\mathbb{E}}[B_t^{2n+1}]=u_n(t,0)=k_n\,t^{n+\frac{1}{2}}$,
where $k_n$ is determined by (\ref{e8}). In particular,
$\hat{\mathbb{E}}[B_t^{3}]=k\,t^{\frac{3}{2}}$, where $k$ satisfies
the following equations:
\begin{equation}
\left \{
\begin{array}{l}
\displaystyle 1
=c^2-c^3\exp(\frac{c^2}{2})\int_c^{\infty}\exp(-\frac{t^2}{2})\,dt
\\[2mm]
\displaystyle
k=-2c^3\exp(\frac{c^2}{2})\\[2mm]
\displaystyle c<0
\end{array}
\right.
\end{equation}
In fact, for each fixed $\sigma\in [0,1]$, we can get
$\hat{\mathbb{E}}[B_t^{3}+lB_t^2]$ for each $l\in\mathbb{R}$ by the
corresponding solution.
\end{remark}

\begin{remark}
Our method can be extended to the functions
\[
\phi(x)=l_1(x^+)^{\alpha}+l_2(x^-)^{\alpha},
\]
where $\alpha >0$ is any given number, $l_1,\ l_2\in\mathbb{R}$. In
fact, this is the only type of functions satisfying $\phi(\lambda
x)= \lambda^\alpha\phi(x)$, $\forall\, \lambda\geq0, x\in
\mathbb{R}$. It is easy to obtain the related Lemma 6, then we get a
method to solve $G$-heat equation (\ref{e1}) with this kind of
initial conditions. Specially, for each $\sigma\in [0,1]$, let
$u_n^{\sigma}$ denote the solution of $G$-heat equation (\ref{e3}),
then $\hat{u}_n^{\sigma}(t,x):=u_n^{\sigma}(t,-x)$ is the solution
of $G$-heat equation (\ref{e1}) with the initial condition
$\phi(x)=-x^{2n+1}$.
\end{remark}

\section{Application to Mathematical Finance}

Let $(\Omega,\mathcal {F},P)$ be a probability space and
$(B_t)_{t\geq0}$ be a $1$-dimensional standard Brownian motion in
this space. We denote by $(\mathcal {F}_t)_{t\geq 0}$ the natural
filtration generated by $B$, i.e.,
\[
\mathcal {F}_t:=\sigma\{ B_s:s\leq t \}\vee \mathcal {N},
\]
where $\mathcal {N}$ denotes the set of all $P$-null subsets in
$\mathcal {F}$. For given $\sigma\in [0,1]$, we denote by $\mathcal
{U}_{\sigma}$ the space of all $\mathcal {F}_t$-adapted controls
with values in $[\sigma, 1]$. For fixed $T>0$, let $(S_t)_{t\leq T}$
be the price of a stock and let $(\log S_t)_{t\leq T}$ satisfy the
following stochastic differential equation
\[
d \log S_t =\mu dt + \sigma_t dB_t,\quad S_0 =1,
\]
where $\mu$ is a constant, $\sigma_.\in \mathcal {U}_{\sigma}$. In
mathematical finance, we often need to calculate (see \cite{Peng-1})
\[
\sup_{\sigma_.\in\mathcal {U}_{\sigma}}E_{P}[\varphi(S_T)].
\]
Specially, for $\varphi (S_T)=(\log S_T)^n$, $n\geq 1$, according to
Remark 4,
\[
\sup_{\sigma_.\in\mathcal {U}_{\sigma}}E_{P}[(\log
S_T)^n]=\hat{\mathbb{E}}[(\mu T +B_T)^n].
\]
Therefore, by Theorem 5, we can get $\sup_{\sigma_.\in\mathcal
{U}_{\sigma}}E_{P}[(\log S_T)^n]$.

\  \  \

{\noindent\bf Acknowledgement}

\  \  \

I would like to thank Professor S.Peng for valuable suggestions and
helpful discussions.

\end{document}